\documentclass[10pt,twoside,a4paper,reqno]{amsart}
\usepackage{amscd,amsmath,amsthm,amsfonts,latexsym,amssymb}

\theoremstyle{plain}
\newtheorem{theo+}           {Theorem}      
\newtheorem{prop+}           {Proposition}  
\newtheorem{coro+}           {Corollary}    
\newtheorem{lemm+}           {Lemma}        
\newtheorem{conjecture}      {Conjecture}

\theoremstyle{definition}
\newtheorem{defi+}           {Definition}   
\newtheorem{prob+}           {Problem}
\newtheorem{not+}            {Notation}

\theoremstyle{remark}
\newtheorem{rema+}           {Remark}       

\newcommand{\bC}{\mathbf C}
\newcommand{\bR}{\mathbf R}
\newcommand{\ve}{\varepsilon}
\newcommand{\De}{\Delta}
\newcommand{\al}{\alpha}

\newcommand{\ze}{\zeta}
\newcommand{\la}{\lambda}
\newcommand{\ka}{\kappa}
\newcommand{\bx}{\mathbf x}
\newcommand{\by}{\mathbf y}
\newcommand{\bbe}{\mathbf e}
\newcommand{\bo}{\mathbf 0}

\newenvironment{theorem}{\begin{theo+}}{\end{theo+}}
\newenvironment{proposition}{\begin{prop+}}{\end{prop+}}
\newenvironment{corollary}{\begin{coro+}}{\end{coro+}}
\newenvironment{lemma}{\begin{lemm+}}{\end{lemm+}}
\newenvironment{remark}{\begin{rema+}}{\end{rema+}}
\newenvironment{definition}{\begin{defi+}}{\end{defi+}}
\newenvironment{notation}{\begin{not+}}{\end{not+}}
\newenvironment{problem}{\begin{prob+}}{\end{prob+}}

\begin{document}
\numberwithin{equation}{section}

\title[Maximal and inextensible polynomials]{Maximal and linearly inextensible 
polynomials}
\author[J.~Borcea]{Julius Borcea}
\address{Department of Mathematics, Stockholm University, SE-106 91 
Stockholm, Sweden}
\email{julius@math.su.se}

\begin{abstract}
Let $S(n,0)$ be the set of monic complex polynomials of degree 
$n\ge 2$ having all their zeros in the closed unit disk and vanishing at 0.
For $p\in S(n,0)$ denote by $|p|_{0}$ the distance from the 
origin to the zero set of $p'$. We determine all $0$-maximal polynomials of 
degree $n$, that is, all polynomials $p\in S(n,0)$ such that 
$|p|_{0}\ge |q|_{0}$ for any 
$q\in S(n,0)$. Using a second order variational method we then 
show that although some of these polynomials are linearly inextensible, 
they are not locally maximal for Sendov's conjecture. 
\end{abstract}

\maketitle

\section{Introduction}\label{s-in}

Let $S_{n}$ be the set of all monic complex polynomials of degree 
$n\ge 2$ having all their zeros in the closed unit disk $\bar{D}$. If 
$p\in S_{n}$ and 
$a\in Z(p)$ then the Gauss-Lucas theorem implies that 
$(a+2\bar{D})\cap Z(p')\neq \emptyset$, where $Z(p)$ and $Z(p')$ denote the 
zero sets of $p$ and $p'$, respectively. In 1959 Sendov conjectured that 
this result may be substantially improved in the following way:
\begin{conjecture}\label{send1}
If $p\in S_{n}$ and $a\in Z(p)$ then $(a+\bar{D})\cap Z(p')\neq \emptyset$.
\end{conjecture}
Sendov's conjecture is widely regarded as one of the main challenges in the 
analytic theory of polynomials. Numerous attempts to verify this conjecture 
have led to over 80 papers, but have met with limited success. We refer to 
\cite{RS}, \cite{Se} and \cite{Sh} for surveys of the results on Sendov's 
conjecture and related questions.

The set $P_n$ of monic complex polynomials of degree $n$ may be 
viewed as a metric space by identifying it with the quotient of 
$\mathbf{C}^n$ by the action of the symmetric group on $n$ elements 
$\Sigma_n$. Indeed, let $\tau:\mathbf{C}^n\rightarrow \mathbf{C}^n/\Sigma_n$ 
denote the orbit map. Let further $p(z)=\prod_{i=1}^{n}(z-z_i)\in P_n$, 
$q(z)=\prod_{i=1}^{n}(z-\zeta_i)\in P_n$ and set
$$\Delta(p,q)=
\min_{\sigma\in \Sigma_n}\max_{1\le i\le n}|z_i-\zeta_{\sigma(i)}|.$$
Then $\Delta$ is a distance function on $P_n$ that 
induces a structure of compact metric space on the set 
$S_n=\{p\in P_n:\Delta(p,z^n)\le 1\}=\tau(\bar{D}^n)$. 
Conjecture~\ref{send1} is therefore an extremum problem 
in the closed unit ball in $P_n$ for the function $d$ given by
$$d:S_{n}\rightarrow [0,2],
\quad p\mapsto d(p)=\max_{z\in Z(p)}\min_{w\in Z(p')}|z-w|.$$
Note that $d(p)$ is the same as the so-called directed (or oriented) 
Hausdorff distance from $Z(p)$ to $Z(p')$ (cf.~\cite{Se}). Since $d$ is 
obviously a continuous function it follows by compactness that 
there exists $p\in S_{n}$ such that $d(p)=\sup_{q\in S_{n}}d(q)$. A 
polynomial with this property is called {\em extremal} for Sendov's 
conjecture. In 1972 Phelps and Rodriguez proposed the following strengthened 
form of Sendov's conjecture (cf.~\cite{PR}):
\begin{conjecture}\label{send2}
If $p\in S_{n}$ is extremal for Sendov's conjecture then 
$p(z)=z^n+e^{i\theta}$ for some $\theta \in \mathbf{R}$.
\end{conjecture}

Variational 
methods are a natural way of dealing with 
Sendov's conjecture. The relatively few properties known to 
hold for locally 
maximal or indeed even extremal polynomials were all deduced by using such 
methods (see \cite{B}, \cite{M1}, \cite{M2}, \cite{SSz}). In this spirit,  
Miller proposed a slightly more general extremum problem in \cite{M1} and 
\cite{M2}. Let $\beta\in \bar{D}$ and denote by $S(n,\beta)$ the 
set of all polynomials in $S_{n}$ having at least one zero at $\beta$. For 
$\alpha\in \mathbf{C}$ and $p\in S(n,\beta)$ let 
$$|p|_{\alpha}=\min_{w\in Z(p')}|\alpha-w|$$ 
and define the $\alpha$-{\em critical circle} to be the circle with center 
$\alpha$ and radius $|p|_{\alpha}$. If $p\in S(n,\beta)$ is such that 
$|p|_{\alpha}\ge |q|_{\alpha}$ for 
any $q\in S(n,\beta)$ then $p$ is said to be {\em maximal with respect to} 
$\alpha$ in $S(n,\beta)$. Maximal polynomials 
with respect to $\alpha$ in 
$S(n,\alpha)$, $\alpha\in \bar{D}$, will be called $\alpha$-{\em maximal} 
for short. A compactness argument again 
shows that maximal polynomials do exist for any $\alpha \in \mathbf{C}$ and 
$\beta\in \bar{D}$ (\cite[Proposition 2.3]{M1}). In 1984 Miller made the 
following conjecture:

\begin{conjecture}\label{miller}
If $p\in S(n,\beta)$ is maximal with respect to $\alpha$ then $p'$ has all 
its zeros on the $\alpha$-critical circle and $p$ has as many zeros as 
possible on the unit circle.
\end{conjecture}

For results pertaining to Miller's conjecture we refer to \cite{B}, 
\cite{M1}, \cite{M2} and \cite{SSz}. 
The relevance of $\alpha$-maximal polynomials in this 
context is quite clear. As shown in \cite[Proposition 2.4]{M1}, if 
$p\in S_{n}$ is an extremal polynomial for Sendov's conjecture then there 
exists $\alpha\in Z(p)$ such that $p$ is $\alpha$-maximal and 
$|p|_{\alpha}=d(p)$. 
Note that {\em \`a priori} there may exist $\alpha\in \bar{D}$ such that if 
$p$ is an 
$\alpha$-maximal polynomial then $|p|_{\alpha}<d(p)$. The $\alpha$-maximal 
polynomials that satisfy $|p|_{\alpha}=d(p)$ are 
particularly interesting since all extremal polynomials are 
necessarily of this type and they may also provide potential 
candidates for local maxima for the function $d$. 
It has been known for some time now that if 
$|\alpha|=1$ then $z^n-\alpha^n$ is the only $\alpha$-maximal polynomial 
(cf.~\cite{Ru}). Moreover, this 
polynomial was shown to be a local maximum for $d$ (\cite{M3}, 
\cite{VZ}). These are actually all the examples of $\alpha$-maximal 
polynomials known so far, 
as no such polynomials were found explicitly for $|\alpha|<1$. 

The main purpose of this paper is to construct and study new examples of 
$\alpha$-maximal polynomials. In \S \ref{s-2} 
we determine all $0$-maximal polynomials 
(Theorem~\ref{class0max}) and investigate some of their properties. 
It turns out 
that all polynomials of this type satisfy Miller's conjecture as well as the 
relation
$|p|_{0}=d(p)$ (Remark~\ref{remmax1} and Corollary~\ref{cormax}). Furthermore,
we show that these polynomials are locally extremal for Sendov's problem in 
$S(n,0)$, i.e., they are local maxima for the restriction of 
the function $d$ to $S(n,0)$ (Proposition~\ref{locmax0}). We conjecture that 
these are in fact all extremal polynomials for Sendov's conjecture 
in the class $S(n,0)$ (Conjecture~\ref{max0}). 

In \S \ref{s-1} we define and characterize a weaker version of local 
maximality, namely the notion of {\em linearly inextensible} polynomial 
(Definition~\ref{lin-in} and Theorems~\ref{nec-suff}--\ref{gen-crit}). 
We then investigate the properties of $0$--maximal polynomials that are  
rotations of the polynomial $p(z)=z^n+z$. If $n\ge 4$ then 
$p$ is locally maximal for a large class of variations of its zeros 
(Proposition~\ref{propvar4}). Moreover, $p$ is locally maximal for the 
restriction of the function $d$ to $S(n,0)$ (Proposition~\ref{locmax0}),  
it is linearly inextensible (Theorem~\ref{thm15}) and satisfies 
Conjecture~\ref{miller} (Remark~\ref{remmax1}). These properties 
and the symmetrical distribution of the 
zeros and critical points of $p$ seem to suggest that if $n\ge 4$ then $p$ 
and its rotations could in fact be locally maximal for Sendov's conjecture in 
the whole class $S_n$.
The discussion in \S \ref{s-1} shows that first order approximations 
are not enough for deciding whether this is true or not. In \S \ref{s-3} 
we use a second order variational method to prove that despite all appearances 
$p$ is not locally maximal for Sendov's conjecture. Indeed,  
Theorems~\ref{secord4} and~\ref{secord5} show that if $n=4$ or $5$ then 
$p$ is actually a kind of inflection point for the function $d$. 
We conjecture that the same is true for all higher degrees 
(Remark~\ref{gencase}). These
results complement those obtained in \S \ref{s-2} and \S \ref{s-1} 
and show quite clearly that 
first order variational methods cannot provide successful ways of dealing with 
Sendov's conjecture in its full generality. 

\subsection*{Acknowledgements} 
The author wishes to thank A.~Meurman and B.~Shapiro for 
stimulating discussions and the anonymous referee for useful suggestions.

\section{Classification of $0$-maximal polynomials}\label{s-2}

In this section we determine all 
$0$-maximal polynomials (Theorem~\ref{class0max}) and study some of 
their properties. 
We show that these polynomials are local maxima for the function $d|_{S(n,0)}$ 
(Proposition~\ref{locmax0}) and we 
conjecture that they are in fact all extremal polynomials for Sendov's 
problem in the class $S(n,0)$ (Conjecture~\ref{max0}). 

We start with a few preliminary results, the first of which 
is well known from the theory of self-inversive polynomials:

\begin{lemma}\label{l-inv}
Let $p(z)=\sum_{k=0}^{n}a_{k}z^k$ be a complex polynomial of degree $n\ge 1$. 
If all the zeros of $p$ lie on the unit circle then $a_{k}\overline{a}_{0}=
a_{n}\overline{a}_{n-k}$, $0\le k\le n-1.$
\end{lemma}

From Lemma~\ref{l-inv} we deduce the following result.

\begin{lemma}\label{l-l2}
Let $q$ be a monic complex polynomial of degree $n\ge 2$ and 
$\alpha\in Z(q)$. If
$R>0$ is such that $|w-\alpha|=R$ for 
any $w\in Z(q')$ then
\begin{equation*}
nR^{2(n-1)}q'(\alpha+z)=q'(\alpha)z^{n-1}\overline{q'(\alpha+R^2\overline{z}^
{\,-1})}
\end{equation*}
for all $z\in \mathbf{C}$, that is,
\begin{equation}\label{lem2}
(n-k-1)!\,nR^{2k}q^{(k+1)}(\alpha)=k!\,q'(\alpha)\,\overline{q^{(n-k)}
(\alpha)},\quad 0\le k\le n-1.
\end{equation}
\end{lemma}

\begin{proof}
By assumption, the polynomial
$$q'(\alpha+Rz)=\sum_{k=0}^{n-1}\frac{R^kq^{(k+1)}(\alpha)}{k!}z^k$$
is of degree at least one and has all its zeros on the unit circle. The 
identities in \eqref{lem2} follow by applying Lemma~\ref{l-inv} to the 
polynomial $q'(\alpha+Rz)$.
\end{proof}

We shall also need the following lemma.

\begin{lemma}\label{l-l3}
If $n$ is an integer greater than two then 
\begin{equation*}
\left\{x\in [1,n-2]\,:\, n^{\frac{2x}{n-1}}(n-x)-n(x+1)=0\right\}=
\left\{\frac{n-1}{2}\right\}.
\end{equation*}
\end{lemma}
\begin{proof}
The assertion is trivially true for $n=3$ and we may therefore assume that 
$n\ge 4$. Note that
$$(-1,n)\ni x\mapsto f(x)=\left(\frac{2x}{n-1}-1\right)\!\log n +\log\!\left(
\frac{n-x}{x+1}\right)$$
is a continuously differentiable function which satisfies 
$$f(0)=f\!\left(\frac{n-1}{2}\right)=f(n-1)=0.$$
It is easily 
seen that the lemma is in fact equivalent to the following statement:
\begin{equation}\label{lem3-1}
f(x)\neq 0 \,\text{ for any }\, x\in [1,n-2]\setminus \left\{\frac{n-1}{2}
\right\}.
\end{equation}
To prove~\eqref{lem3-1} note that the equation $f'(x)=0$ has at most two real 
roots in the interval $[0,n-1]$ since
$$f'(x)=\frac{2(\log n) x^2-2[(n-1)\log n] x+n^2-2n\log n -1}
{(n-1)(x+1)(x-n)}.$$
Thus by Rolle's theorem the equation $f(x)=0$ cannot have more than the three
roots listed above in the interval $[0,n-1]$, which establishes~\eqref{lem3-1}.
\end{proof}

We are now ready to prove the main result of this section:

\begin{theorem}\label{class0max}
The $0$-maximal polynomials of degree $n\ge 2$ are given by:
\begin{itemize}
\item[(i)] $z^{2m}+e^{i2\theta}z$ if $n=2m$ and $m\ge 1$, where 
$\theta\in \mathbf{R}$.
\item[(ii)] $z^{2m+1}+\lambda e^{i\theta}z^{m+1}+e^{i2\theta}z$ if 
$n=2m+1$ and $m\ge 1$, where $\lambda,\theta\in \mathbf{R}$ and 
$|\lambda|\le \frac{2\sqrt{2m+1}}{m+1}$.
\end{itemize}
\end{theorem}
\begin{proof}
Let $p$ be a $0$-maximal polynomial. Since $p\in S(n,0)$, we may write
\begin{eqnarray*}
&&p(z)=z^n+\sum_{k=1}^{n-1}a_{k}z^k=z\prod_{i=1}^{n-1}(z-z_{i}),
\text{ where } |z_{i}|\le 
1,\,1\le i\le n-1,\\
&&p'(z)=nz^{n-1}+\sum_{k=1}^{n-1}ka_{k}z^{k-1}=n\prod_{j=1}^{n-1}(z-w_{j}).
\end{eqnarray*}
By comparing $p$ with the polynomial $q(z):=z^{n}-z\in S(n,0)$ we get
$$\min_{1\le j\le n-1}|w_{j}|=|p|_{0}\ge |q|_{0}=\left(\!\frac{1}{n}\!\right)^
{\frac{1}{n-1}},$$
which combined with the identity $\prod_{i=1}^{n-1}z_{i}=(-1)^{n-1}p'(0)=
n\prod_{j=1}^{n-1}w_{j}$ yields 
$\prod_{i=1}^{n-1}|z_{i}|=n\prod_{j=1}^{n-1}|w_{j}|\ge 1$. Since $|z_{i}|
\le 1$ for $1\le i\le n-1$ we deduce that
\begin{equation}\label{thm2-1}
|z_{i}|=1,\,1\le i\le n-1,\,\text{ and }\,|w_{j}|=\left(\!\frac{1}{n}\!
\right)^{\frac{1}{n-1}},\,1\le j\le n-1.
\end{equation}
It follows in particular that $|a_{1}|=|\prod_{i=1}^{n-1}z_{i}|=1$. Hence
\begin{equation}\label{thm2-2}
a_{1}=e^{i2\theta}\text{ for some }\theta\in \mathbf{R}.
\end{equation}
From \eqref{thm2-1} and Lemma~\ref{l-inv} applied to the polynomial 
$z^{-1}p(z)$ we get
$$\overline{a}_{n-k}=e^{-i2\theta}a_{k+1},\quad 1\le k\le n-2,$$
while \eqref{thm2-1} and Lemma~\ref{l-l2} applied to the polynomial $p$ with 
$\alpha=0$ and $R=\left(\!\frac{1}{n}\!\right)^{\frac{1}{n-1}}$ 
imply that
$$\overline{a}_{n-k}=\frac{n(k+1)}{n-k}\left(\!\frac{1}{n}\!\right)^{
\frac{2k}{n-1}}e^{-i2\theta}a_{k+1},\quad 1\le k\le n-2.$$
Thus
\begin{equation}\label{thm2-3}
\left[n^{\frac{2k}{n-1}}(n-k)-n(k+1)\right]a_{k+1}=0,\quad 1\le k\le n-2,
\end{equation}
and then by Lemma~\ref{l-l3} we get that $a_{k}=0$ for $2\le k\le n-1$ if $n$ 
is even. 

Let now $n=2m+1,\,m\ge 1$. In this case Lemma~\ref{l-l3} and \eqref{thm2-3} 
imply that $a_{k}=0$ for $2\le k\le n-1,\,k\neq m+1$, so that
\begin{equation}\label{thm2-4}
p(z)=z^{2m+1}+a_{m+1}z^{m+1}+e^{i2\theta}z,
\end{equation}
where $\theta$ is as in \eqref{thm2-2}. From \eqref{thm2-1} and 
\eqref{thm2-4} we obtain
$$\sum_{i=1}^{2m}z_{i}^m+2ma_{m+1}+e^{i2\theta}\sum_{i=1}^{2m}
\overline{z}_{i}^m=0.$$
Note that by Newton's identities one has actually that
$$\sum_{i=1}^{2m}z_{i}^m=-ma_{m+1}-\sum_{k=1}^{m-1}a_{m+k+1}\bigg(\!
\sum_{j=1}^{2m}z_{j}^{k}\!\bigg)=-ma_{m+1}.$$
We deduce from the last two formulas that $a_{m+1}=e^{i2\theta}
\overline{a}_{m+1}$, which together with \eqref{thm2-4} implies that the 
polynomial $p$ must be of the form
$$p(z)=z^{2m+1}+\lambda e^{i\theta}z^{m+1}+e^{i2\theta}z,\,\text{ where }\,
\lambda\in \mathbf{R}.$$
An elementary computation shows that if $\theta,s,t\in \mathbf{R}$ with $t>0$ 
then the roots of the equation $x^2+se^{i\theta}x+te^{i2\theta}=0$ 
cannot have the same absolute value unless $s^2\le 4t$. By \eqref{thm2-1} all 
the zeros of the polynomial
$$p'(z)=(2m+1)z^{2m}+(m+1)\lambda e^{i\theta}z^{m}+e^{i2\theta}$$
have the same absolute value. Now using the substitution $x=z^m$ and the 
above-mentioned result on the roots of second degree equations we see that 
this cannot happen unless $|\lambda|\le \frac{2\sqrt{2m+1}}{m+1}$.
To summarize, we have shown that if $p\in S(n,0)$ then $|p|_{0}<\left(\!
\frac{1}{n}\!\right)^{\frac{1}{n-1}}$ unless $p$ is of the form (i) or 
(ii) according to the parity of $n$. On the other hand, it is easily checked 
that a polynomial of the form (i) or (ii) satisfies 
$|p|_{0}=\left(\!\frac{1}{n}\!\right)^{\frac{1}{n-1}}$, which completes the 
proof of the theorem.
\end{proof}

\begin{remark}\label{rem-newref}
The estimate $|p|_0\le \left(\!\frac{1}{n}\!\right)^{\frac{1}{n-1}}$ for  
$p\in S(n,0)$ appears also in \cite[Lemma~4]{BRS} but the 
cases of equality are not dealt with in {\em loc.~cit.}
\end{remark}

\begin{remark}\label{remmax1}
Let $p$ be a $0$-maximal polynomial of degree $n\ge 2$. By 
Theorem~\ref{class0max} all the 
critical points of $p$ lie on the $0$-critical circle $|z|=\left(\!
\frac{1}{n}\!\right)^{\frac{1}{n-1}}$ and all the zeros of $p$ except $0$ lie 
on the unit circle. Thus, all $0$-maximal polynomials satisfy Miller's 
conjecture (Conjecture~\ref{miller}). Note also that if $p$ is $0$-maximal 
then it has only simple zeros. The same is true for $p'$ except 
when $n$ is odd and $p$ is of the form (ii) with 
$|\lambda|=\frac{4\sqrt{n}}{n+1}$, in which case $p'$ has $\frac{n-1}{2}$ 
distinct zeros each of multiplicity two.
\end{remark}

\begin{corollary}\label{cormax}
If $p$ is a $0$-maximal polynomial of degree $n\ge 2$ then $d(p)=|p|_0$.
\end{corollary}

\begin{proof}
Let $z\in Z(p)\setminus \{0\}$, so that $|z|=1$ by Remark~\ref{remmax1}. 
Denote by $H_z$ the closed half-plane that contains $z$ and is bounded by 
the line passing through $\frac{z}{2}$ which is orthogonal to the 
segment $[0,z]$. Since $p(z)=0=p(0)$ it follows from a well-known 
consequence of the Grace-Heawood theorem that there exists 
$w\in Z(p')\cap H_z$ (see, e.g., the supplement to 
Theorem 4.3.1 in \cite{RS}). By 
Remark~\ref{remmax1} again one has that $|w|=|p|_0=\left(\!
\frac{1}{n}\!\right)^{\frac{1}{n-1}}\ge \frac{1}{2}$. It is then geometrically 
clear that $|z-w|\le |p|_0$ and thus 
$\min_{\omega\in Z(p')}|z-\omega|\le |p|_0$, which proves the desired result. 
\end{proof}  

Theorem~\ref{class0max} and Corollary~\ref{cormax} do not automatically imply 
that $0$-maximal polynomials are global maxima for the restriction of the 
function $d$ to the subset $S(n,0)$ of $S_n$. For this one would have to 
show that $d(p)\le \left(\!\frac{1}{n}\!\right)^{\frac{1}{n-1}}$ 
whenever $p\in S(n,0)$, which we actually conjecture to be true: 

\begin{conjecture}\label{max0}
If $p\in S(n,0)$ then $d(p)\le \left(\!\frac{1}{n}\!\right)^{\frac{1}{n-1}}$ 
and equality occurs if and only if $p$ is as in Theorem~\ref{class0max}.
\end{conjecture}

In other words, Conjecture~\ref{max0} asserts that Theorem~\ref{class0max} 
gives all extremal polynomials for Sendov's problem restricted to the class 
$S(n,0)$. We should point out that the best estimate known so far for 
$p\in S(n,0)$ is $d(p)<1$ (see~\cite[Theorem 7.3.6]{RS}). 
Our next result 
shows that the polynomials in Theorem~\ref{class0max} are locally 
maximal for Sendov's problem in the class $S(n,0)$, so that 
Conjecture~\ref{max0} is at least locally valid.

\begin{proposition}\label{locmax0}
If $p$ is a polynomial as 
in Theorem~\ref{class0max} there exists $\varepsilon>0$ such 
that whenever $q\in S(n,0)$ and $\Delta(q,p)<\varepsilon$ one has 
$d(q)\le \left(\!\frac{1}{n}\!\right)^{\frac{1}{n-1}}$. Moreover, the 
latter inequality is 
strict unless $q$ is itself a polynomial as in Theorem~\ref{class0max}.
\end{proposition}

\begin{proof}
The statement is trivially true if $n=2$ and so we may assume that $n\ge 3$. 
If $p$ is a polynomial as in Theorem~\ref{class0max} then the proof of 
Corollary~\ref{cormax} shows that 
$|p|_{z}<|p|_0=d(p)$ for $z\in Z(p)\setminus \{0\}$. Hence there 
exists $\varepsilon=\varepsilon(p)>0$ such that $|q|_{\zeta}<|q|_0=d(q)$ for 
$\zeta\in Z(q)\setminus \{0\}$ whenever $q\in S(n,0)$ and 
$\Delta(q,p)<\varepsilon$. The result is now a consequence of 
Theorem~\ref{class0max}.
\end{proof}

\begin{remark}\label{remmax2}
It is not difficult to see that that $\max_{p\in S(n,1)}|p|_{1}=1$ 
(cf.~\cite{Ru}). On the other hand, Theorem~\ref{class0max} shows that 
$\max_{p\in S(n,0)}|p|_{0}$ 
increases to $1$ as $n\rightarrow \infty$. It would be interesting to know 
whether there exist $\alpha\in D$ 
and $c\in (0,1)$ such that $\max_{p\in S(n,\alpha)}|p|_{\alpha}\le c$ for any 
$n\ge 2$.
\end{remark}

To end this section let us point out that Theorem~\ref{class0max} solves 
in fact the following more general extremum problem 
concerning the distribution of zeros and critical points of complex 
polynomials:

\begin{problem}\label{pb1}
Let $n\ge 2$, $a\in \mathbf{C}$, and $R>0$. Find the largest constant 
$\rho:=\rho(a,n,R)$ with the property that for any complex polynomial $p$ of 
degree $n$ satisfying $p(a)=0$ and $\min_{w\in Z(p')}|w-a|\ge R$ one has 
$\max_{z\in Z(p)}|z-a|\ge \rho$.
\end{problem}

It is easy to see that $\rho$ is invariant under translations in the complex 
plane so that it does not depend on $a$. Moreover, the Gauss-Lucas theorem 
clearly implies that $\rho >R$. 
On the other hand, by considering $p(z)=(z-a)^n+nR^{n-1}(z-a)$ we see that 
$\rho \le Rn^{\frac{1}{n-1}}$. Essentially the same computations as in the 
proof of Theorem~\ref{class0max} show that in fact 
$\rho=Rn^{\frac{1}{n-1}}$ and the extremal cases are as follows.

\begin{theorem}
In the above notations one has $\rho=Rn^{\frac{1}{n-1}}$ and this 
value is attained only for the following polynomials:
\begin{itemize}
\item[(i)] $(z-a)^{2m}+2mR^{2m-1}e^{i\theta}(z-a)$ if $n=2m$ and $m\ge 1$,
where $\theta\in \mathbf{R}$.
\item[(ii)] $(z-a)^{2m+1}+\lambda \sqrt{2m+1}R^{m}e^{i\theta}(z-a)^{m+1}+
(2m+1)R^{2m}e^{i2\theta}(z-a)$ if $n=2m+1$ and $m\ge 1$, where 
$\lambda,\theta\in \mathbf{R}$ and $|\lambda|\le \frac{2\sqrt{2m+1}}{m+1}$.
\end{itemize}
\end{theorem}

\section{First order expansions and linearly inextensible 
polynomials}\label{s-1}

In the previous section we have already mentioned the notion of local maximum 
for the Sendov problem. A formal definition of this notion is as follows.

\begin{definition}\label{locmaxdef}
Given $\ve>0$ and $p\in S_n$ we define the $\ve$-neighborhood of $p$ 
to be 
$V_{\ve}(p)=\{q\in S_n\,:\,\De(q,p)<\ve\}.$ A polynomial $p\in S_n$ is 
called {\em locally maximal for Sendov's conjecture} if 
it is a local maximum for the function $d$, i.e., there exists 
$\varepsilon>0$ such that $d(q)\le d(p)$ for any $q\in V_{\ve}(p)$.
\end{definition}

We shall now define and study a weaker version of this notion.

\begin{definition}\label{lin-in}
A polynomial $p\in S_n$ is {\em linearly inextensible} if there 
exist $\ve>0$ and constants $C\ge 0$, $\ka>0$ depending only on 
$p$ such that $d(q)\le d(p)+C\De(q,p)^{1+\ka}$ for all $q\in V_{\ve}(p)$. 
Otherwise $p$ is said to be {\em linearly extensible}.
\end{definition}

Clearly, a locally maximal polynomial is necessarily linearly inextensible 
while a linearly extensible polynomial cannot be locally maximal. Necessary 
criteria for local maximality in certain generic cases were established in 
\cite{B} and \cite{M1}. 

Definition~\ref{lin-in} is 
motivated mainly by two facts. On the one hand, a close examina\-tion of the 
aforementioned criteria shows that they are in fact both 
necessary and sufficient conditions for linear inextensibility in the generic 
cases under consideration (Theorems~\ref{nec-suff}--\ref{gen-crit}). On the 
other hand, linear inextensibility can be decided 
by means of first order expansions while local maximality may require 
higher order approximations. 
The examples below show that these are actually 
crucial differences between the notions in Definitions~\ref{locmaxdef} 
and~\ref{lin-in} and that first order variational methods 
are clearly not enough for studying local maxima for Sendov's conjecture.

We need some further preliminaries. Let us first recall 
\cite[Definition 2.19]{M1}:

\begin{definition}\label{def13}
A complex $m\times n$ matrix $M=(m_{ij})$ is {\em positively 
singular} if there exist $\mu_{1},\ldots,\mu_{m}\ge 0$, not all $0$, so that 
$\sum_{i=1}^{m}\mu_{i}m_{ij}=0$ for $1\le j\le n$.
\end{definition}

Given a vector $\by=(y_1,\ldots,y_n)\in\bC^n$ we shall formally write 
$\Re(\by)>0$ whenever the inequalities $\Re(y_i)>0$, $1\le i\le n$, are 
satisfied.

\begin{theorem}\label{thm14}
Let $M$ be a complex $m\times n$ matrix and let $\bx\in 
\mathbf{C}^n$ denote a vector of complex unknowns. The system 
$\Re\left(M\bx\right)>0$ is solvable 
if and only if $M$ is not positively singular.
\end{theorem} 

\begin{remark}
Theorem~\ref{thm14} is a generalization of what is usually called the 
fundamental theorem of linear programming, see e.g.~\cite[Theorem 22.2]{Ro} 
and \cite[Lemma 2.20]{M1}.
\end{remark}

\begin{notation}\label{not-1}
For $p\in S_{n}$ we shall use the following generic notations:
\begin{equation}\label{poly1}
\begin{split}
&p(z)=\prod_{i=1}^{n}(z-z_{i}),\quad a=z_{1},\quad |p|_{a}=d(p),\quad p'(z)=
n\prod_{j=1}^{n-1}(z-w_{j}),\\
&r=r(a)\in\{1,\ldots,n-1\}\text{ and }1\le s_1<\ldots<s_m\le n\text{ are such 
that}\\
&|w_{j}-a|=|p|_{a}\text{ for }1\le j\le r,\quad |w_{j}-a|>|p|_{a}\text{ for }
j\ge r+1,\\
&|z_{s_k}|=1\text{ for }1\le k\le m,\quad |z_i|<1\text{ for }
i\in\{1,\ldots,n\}\setminus\{s_1,\ldots,s_m\}.
\end{split}
\end{equation}
\end{notation}

The following result combines Lemmas 2.1 and 2.3 in \cite{B}.

\begin{lemma}\label{lem19}
If $p$ is as in \eqref{poly1} and has only simple zeros 
there exist neighborhoods $U, V\subset \mathbf{C}^n$ of the points 
$u=(a,w_{1},\ldots,w_{n-1})$ and 
$v=(a,z_{2},\ldots,z_{n})$, respectively, such that
$$U\ni (\alpha,\omega_{1},\ldots,\omega_{n-1})\mapsto (\zeta_{1},\zeta_{2},
\ldots,\zeta_{n})\in V$$
is an analytic function, where $\zeta_{2},\ldots,\zeta_{n}$ are the (simple) 
zeros different from $\alpha$ of the polynomial 
$n\!\int_{\alpha}^{z}\prod_{j=1}^{n-1}(w-\omega_{j})dw$ and
$\zeta_{1}=\zeta_{1}(\alpha,\omega_{1},\ldots,\omega_{n-1})\equiv\alpha$. 
Moreover,  one has
$$\frac{\partial \zeta_{k}}{\partial \omega_{l}}\bigg|_{u}=
\frac{1}{p'(z_{k})}\int_{a}^{z_{k}}\frac{p'(w)}{w-w_{l}}dw,\quad 2\le k\le n, 
1\le l\le n-1.$$
\end{lemma}

From Lemma~\ref{lem19} and the inverse function theorem one immediately gets:

\begin{lemma}\label{invfcn}
Let $p$ be as in \eqref{poly1} and assume that both $p$ and $p'$ have simple 
zeros. Then in the notations of Lemma~\ref{lem19} the functions 
$\omega_{1},\ldots,\omega_{n-1}$ are locally analytic in 
$\zeta_{1},\zeta_{2},\ldots,\zeta_{n}$ and 
$$\frac{\partial \omega_j}{\partial \zeta_i}\bigg|_{v}=
-\frac{p(w_j)}{(w_j-z_i)^2p''(w_j)},\quad 1\le i\le n,\,1\le j\le n-1.$$
\end{lemma}

\begin{remark}\label{byprod}
An interesting byproduct of Lemmas~\ref{lem19} and~\ref{invfcn} is the 
following ``gene\-ralized Cauchy determinant'' property for complex 
polynomials: let $p$ be a complex polynomial 
of degree $n\ge 2$ such that both $p$ and $p'$ have simple zeros, set 
$Z(p)=\{z_1,\ldots,z_n\}$, $Z(p')=\{w_1,\ldots,w_{n-1}\}$ and define an 
$(n-1)\times (n-1)$ matrix $B(p)=((z_i-w_j)^{-2})$, $1\le i,j\le n-1$. Then
$\det(B(p))\neq 0$ (cf.~\cite[Theorem 1.3]{B1}).
\end{remark}

\begin{notation}\label{not-2}
Denote by $||\cdot||$ the norm on $\bC^n$ given by 
$||\bx||=\max_{1\le i\le n}|x_i|$ for $\bx=(x_1,\ldots,x_n)\in \bC^n$. 
Given $\bbe=(\ve_1,\ldots,\ve_n)\in\bC^n$ and a polynomial $p\in S_n$ as 
in \eqref{poly1} define $z_{i}(\bbe)=z_{i}+\ve_i$, $1\le i\le n$, and set 
\begin{equation}\label{poly2}
p_{\bbe}(z)=\prod_{i=1}^{n}(z-z_{i}(\bbe)),\quad p_{\bbe}'(z)=
n\prod_{j=1}^{n-1}(z-w_{j}(\bbe)).
\end{equation}
Note that $z_{i}(\bo)=z_{i}$, $1\le i\le n$, 
$w_{j}(\bo)=w_{j}$, $1\le j\le n-1$, and that if $p$ has simple zeros then 
$p_{\bbe}$ has simple zeros and 
$\De(p_{\bbe},p)=||\bbe||$ for all $\bbe\in\bC^n$ with $||\bbe||$ sufficiently
small. Clearly, $p_{\bbe}\in S_n$ if and only if $|z_i+\ve_i|\le 1$, 
$1\le i\le n$.
\end{notation}

Lemma~\ref{invfcn} and straightforward computations yield the following result.

\begin{proposition}\label{prop11}
Let $p$ be as in \eqref{poly1} and such that both $p$ and $p'$ have simple 
zeros. For any $\bbe=(\ve_1,\ldots,\ve_n)\in\bC^n$ with 
$||\bbe||$ sufficiently small one has
$$|w_{j}(\bbe)-z_{1}(\bbe)|=|p|_{a}\!\left[1+\sum_{i=1}^{n}
\Re\left(\al_{i}(w_{j};a)\ve_{i}\right)
+\mathcal{O}\left(||\bbe||^2\right)\right],\quad 
1\le j\le r(a),$$
where 
\begin{align*}
&\al_{1}(w_{j};a)=-\frac{1}{w_{j}-a}\!\left[1+\frac{p(w_{j})}{(w_{j}-a)^{2}
p''(w_{j})}\right],\\
&\al_{i}(w_{j};a)=-\frac{p(w_{j})}{(w_{j}-a)(w_{j}-z_{i})^{2}p''(w_{j})},
\quad 2\le i\le n.
\end{align*}
\end{proposition}

\begin{notation}\label{not-3}
Given a polynomial $p$ as in \eqref{poly1} such that both $p$ and $p'$ have 
simple zeros we define an $(r(a)+m)\times n$ matrix $A(p;a)=(a_{ij})$ by 
setting
\begin{equation}\label{matr-a}
a_{ij}=
\begin{cases}
\al_{j}(w_{i};a)&\text{ if } 1\le i\le r(a),\,1\le j\le n,\\
-\overline{z}_{s_k}&\text{ if } i=r(a)+k\text{ and }j=s_k,\,1\le k\le m,\\
0&\text{ otherwise}.
\end{cases}
\end{equation}
\end{notation}

We can now prove a necessary and sufficient criterion for linear 
inextensibility.

\begin{theorem}\label{nec-suff}
Let $p\in S_n$ be a polynomial as in \eqref{poly1} with $n\ge 3$. Assume that 
both $p$ and 
$p'$ have simple zeros and that $|p|_{z_i}<d(p)$ for $2\le i\le n$. Then 
$p$ is linearly inextensible if and only if $A(p;a)$ is positively singular.
\end{theorem}

\begin{proof}
Suppose that $p$ is linearly inextensible and that there exists $\bbe\in\bC^n$ 
such that $\Re(A(p;a)\bbe)>0$. Let $t>0$ and construct the polynomial 
$p_{t\bbe}(z)$ as in \eqref{poly2}. Then $p_{t\bbe}$ has simple zeros and 
$\De(p_{t\bbe},p)=t||\bbe||$ for all sufficiently small $t$. By 
Proposition~\ref{prop11} there exists a constant $\la>0$ depending only on 
$p$ (and $\bbe$) such that
$$d(p_{t\bbe})=|p_{t\bbe}|_{z_1(t\bbe)}
=\min_{1\le j\le r(a)}|w_j(t\bbe)-z_1(t\bbe)|>|p|_{a}(1+\la t)$$
for all small $t>0$, which contradicts the fact that $p$ is linearly 
inextensible. Hence the system $\Re(A(p;a)\bx)>0$, $x\in\bC^n$, cannot be 
solvable and so by Theorem~\ref{thm14} $A(p;a)$ must be positively singular.

To establish the sufficiency part note first that if $z,\ve\in\bC$ are 
such that $\ve\neq 0$, $|z|=1$ and $|z+\ve|\le 1$ then 
$\Re(\ve\overline{z})<0$. Recall~\eqref{matr-a}, let 
$1\le k_1<\ldots<k_l\le m$ be an ordered $l$-tuple of indices, where 
$l\le m$, and denote by
$A(p;a)[k_1,\ldots,k_l]$
the $(r(a)+m-l)\times (n-l)$ matrix obtained by deleting rows 
$r(a)+k_1,\ldots,r(a)+k_l$ and columns $s_{k_1},\ldots,s_{k_l}$ from $A(p;a)$.
Suppose now that $A(p;a)$ is positively singular. Then so is any 
matrix of the form $A(p;a)[k_1,\ldots,k_l]$ provided it is not the empty 
matrix, which can only occur if $l=m=n$.

Let $\bbe=(\ve_1,\ldots,\ve_n)\in\bC^n\setminus\{\bo\}$ be such that 
$||\bbe||$ is arbitrarily small and consider the 
polynomial $p_{\bbe}(z)$ as in \eqref{poly2}, on which we impose the condition
$p_{\bbe}\in S_n$. Set 
$$K(\bbe)=\{k_1(\bbe),\ldots,k_l(\bbe)\}
:=\{i\in\{1,\ldots,m\}\,:\,\ve_{s_i}=0\},$$ 
so that $l\le n-1$ since $\bbe\neq \bo$. Assuming as we may that 
$k_1(\bbe)<\ldots<k_l(\bbe)$ let us construct the corresponding (non-empty) 
matrix $A(p;a)[k_1(\bbe),\ldots,k_l(\bbe)]$ as above, where it is understood 
that $A(p;a)[k_1(\bbe),\ldots,k_l(\bbe)]=A(p;a)$ if $K(\bbe)=\emptyset$. Since 
$\Re(\ve_{s_i}\overline{z}_{s_i})<0$ for $i\in \{1,\ldots,m\}\setminus K(\bbe)$
and $A(p;a)[k_1(\bbe),\ldots,k_l(\bbe)]$ is positively singular it follows 
from~\eqref{matr-a} and Theorem~\ref{thm14} that there exists 
$j\in\{1,\ldots,r(a)\}$ such that
$$\sum_{i=1}^{n}\Re(\al_i(w_j;a)\ve_i)\le 0.$$
Proposition~\ref{prop11} implies that there exists a non-negative 
constant $C$ depending only on $p$ such that 
$$d(p_{\bbe})=|p_{\bbe}|_{z_1(\bbe)}=\min_{1\le j\le r(a)}|w_j(\bbe)-z_1(\bbe)|
\le |p|_a\!\left(1+C||\bbe||^2\right),$$
which shows that $p$ is linearly inextensible.
\end{proof}

The following generalization of Theorem~\ref{nec-suff} is obtained in 
similar fashion and we therefore state it without proof.

\begin{theorem}\label{gen-crit}
Let $p\in S_n$ be as in \eqref{poly1} with $n\ge 3$ 
and assume that both $p$ and 
$p'$ have simple zeros. Suppose further that the $z_i$'s are labeled so that 
$|p|_{z_i}=d(p)$ for $1\le i\le t$ and $|p|_{z_i}<d(p)$ for $i\ge t+1$, 
where $t\in\{1,\ldots,n\}$. Then 
$p$ is linearly in\-extensible if and only if 
$A(p;z_i)$, $1\le i\le t$, are positively singular matrices.
\end{theorem}

\begin{remark}\label{other-crit}
Let 
\begin{equation*}
\begin{split}
&\hat{S}_n=\{p\in S_n\,:\,p\text{ and }p'\text{ have simple zeros}\}
\text{ and}\\ 
&\tilde{S}_n=\{p\in S_n\,:\,\text{all zeros of }p\text{ on the unit circle are 
simple}\},\\
\end{split}
\end{equation*}
so that $\hat{S}_n\subsetneq \tilde{S}_n$. The necessity part in 
Theorem~\ref{nec-suff} viewed as a necessary condition
for local maximality in the class $\hat{S}_n$ was established 
in \cite{M1}. A necessary criterion for local maximality in the larger 
class $\tilde{S}_n$ was obtained in \cite{B}. It is not difficult to see that
the latter actually gives necessary and essentially also sufficient 
conditions for linear inextensibility in the class $\tilde{S}_n$.
\end{remark}

The $0$-maximal polynomials given in Theorem~\ref{class0max} provide natural 
examples of linearly inextensible polynomials of arbitrary degree:

\begin{theorem}\label{thm15}
Let $\theta \in \mathbf{R}$,  $n\ge 3$, and $p(z)=z^n+e^{i\theta}z$. Then 
$|p|_{0}=d(p)$, $|p|_{\ze}<d(p)$ for $\ze\in Z(p)\setminus \{0\}$ and 
$p$ is linearly inextensible.
\end{theorem}

\begin{proof}
It is enough to prove the statement for the polynomial $p(z)=z^n-z$ since 
$z^n+e^{i\theta}z=e^{-in\alpha}p(e^{i\alpha}z)$, where 
$\alpha=\frac{\pi-\theta}{n-1}$. Recall the notations of 
\eqref{poly1} and set
$$(a=)\,z_{1}=0,\,z_{j}=e^{\frac{2\pi i(j-2)}{n-1}},\,2\le j\le n,\quad
w_{k}=\left(\!\frac{1}{n}\!\right)^{\!\frac{1}{n-1}}\!
e^{\frac{2\pi i(k-1)}{n-1}},\,1\le k\le n-1.$$
It is easy to check that $|p|_{0}=d(p)$, $|p|_{\ze}<d(p)$ for 
$\ze\in Z(p)\setminus \{0\}$ and $r(0)=n-1$, so that the matrix $A(p;a)$ 
defined in \eqref{matr-a} is actually a $2(n-1)\times n$ matrix. 
Proposition~\ref{prop11} and elementary computations yield
$$\al_{1}(w_{i};a)=-\frac{n-1}{nw_{i}},\quad \al_{j}(w_{i};a)=
\frac{w_{i}}{n(w_{i}-z_{j})^2},\quad 1\le i\le n-1,\,2\le j\le n.$$
Using the fact that for $z\in \mathbf{C}\setminus Z(p')$ one has the 
well-known identities
\begin{eqnarray*}
&&\frac{p''(z)}{p'(z)}=\sum_{i=1}^{n-1}\frac{1}{z-w_{i}}\,\text{ and}\\
&&\frac{p^{(3)}(z)p'(z)-(p''(z))^2}{(p'(z))^2}=\frac{d}{dz}\!
\left[\frac{p''(z)}{p'(z)}\right]=-\sum_{i=1}^{n-1}\frac{1}{(z-w_{i})^2}
\end{eqnarray*}
one can show that if $2\le j\le n$ then
\begin{eqnarray*}
\sum_{i=1}^{n-1}\frac{w_{i}}{(w_{i}-z_{j})^2}&=&\sum_{i=1}^{n-1}
\frac{1}{w_{i}-z_{j}}+z_{j}\sum_{i=1}^{n-1}\frac{1}{(w_{i}-z_{j})^2}\\
&=&-\frac{p''(z_{j})}{p'(z_{j})}+z_{j}\!\left[\frac{(p''(z_{j}))^2-p^{(3)}
(z_{j})p'(z_{j})}{(p'(z_{j}))^2}\right]=\frac{n}{z_{j}}.
\end{eqnarray*}
It follows that
$$\sum_{i=1}^{n-1}\al_1(w_{i};a)
=-\frac{n-1}{n}\sum_{i=1}^{n-1}\frac{1}{w_{i}}=0,\quad  
\sum_{i=1}^{n-1}\al_{j}(w_{i};a)-\bar{z}_{j}
=\frac{1}{z_{j}}-\bar{z}_{j}=0,\quad 2\le j\le n,$$
which shows that $A(p;a)$ is a positively singular matrix. By 
Theorem~\ref{nec-suff} the polynomial $p(z)=z^n-z$ cannot be linearly 
extensible.
\end{proof}

We can now give a concrete example of a linearly inextensible polynomial 
which is not a local maximum for Sendov's conjecture.

\begin{lemma}\label{countlm3}
If $\theta\in \mathbf{R}$ then $p(z)=z^3+e^{i\theta}z$ is 
not a local maximum for $d$.
\end{lemma}

\begin{proof}
Note that $d(p)=|p|_{0}=
\frac{1}{\sqrt{3}}$ and that by a rotation we may 
assume that $p(z)=z^3-z$. Let $t\in [0,1]$ and set
$$q_{t}(z)=(z-it)(z^2-1).$$
Then $q_{t}\in S_{3}$ and elementary computations show that for all small 
$t>0$ one has
$$d(q_{t})=|q_{t}|_{it}=\sqrt{\frac{1+t^2}{3}}>d(p),$$
which proves the lemma.
\end{proof}

\begin{remark}\label{otherex}
In \S \ref{s-3} we shall construct explicit examples of linearly inextensible 
polynomials of degree greater than $3$ which are not local maxima for $d$.
\end{remark}

In the remainder of this section we establish further properties of the 
$0$-maximal and linearly inextensible polynomials of the form $z^n-z$ 
with $n\ge 4$. 

\begin{notation}\label{not-li}
Recall the notations used in the proof of Theorem~\ref{thm15} and set
\begin{equation}\label{c-ex1}
\begin{split}
&p(z)=z^n-z=\prod_{j=1}^{n}(z-z_j)\text{ and }
p'(z)=n\prod_{k=1}^{n-1}(z-w_k),\text{ where }z_{1}=0,\\
&z_{j}=e^{\frac{2\pi i(j-2)}{n-1}},\,2\le j\le n, 
\text{ and } w_{k}=\left(\!\frac{1}{n}\!\right)^{\!\frac{1}{n-1}}\!
e^{\frac{2\pi i(k-1)}{n-1}},\,1\le k\le n-1.
\end{split}
\end{equation}
Let $\kappa$ be a fixed positive number and 
$\bbe=(\varepsilon_1,\ldots,\varepsilon_n)\in \mathbf{C}^n$ be such 
that $|\varepsilon_j|\le |\varepsilon_1|^{1+\kappa}$ for $2\le j\le n$ and let
$p_{\bbe}$ be as in~\eqref{poly2}, that is,
$p_{\bbe}(z)=\prod_{j=1}^{n}(z-(z_j+\varepsilon_j))$.
\end{notation}

Using the above notations we can show the following result:

\begin{proposition}\label{propvar4}
If $n\ge 4$ then for all sufficiently small $|\varepsilon_1|$ one has
$$\min_{\omega\in Z(p_{\bbe}')}|z_1+\varepsilon_1-\omega|\le
\min_{\omega\in Z(p')}|z_1-\omega|
-\left[\cos\!\left(\frac{\pi}{n-1}\right)\right]|\varepsilon_1|.$$
Thus, if $\varepsilon_j=0$ for $2\le j\le n$ or, more generally, if 
$|z_j+\varepsilon_j|\le 1$ for $2\le j\le n$ then for all sufficiently small 
$|\varepsilon_1|>0$ one has $p_{\bbe}\in S_n$ and $d(p_{\bbe})<d(p)$.
\end{proposition}
\begin{proof}
Note that $|w_k|=\left(\!\frac{1}{n}\!\right)^{\!\frac{1}{n-1}}=|p|_{z_1}=
d(p)$, $1\le k\le n-1$, and denote the zeros of $p_{\bbe}'(z)$ by 
$\omega_k=\omega_k(\varepsilon_1,\ldots,\varepsilon_n)$, $1\le k\le n-1$. We 
assume that these are labeled so that $\omega_k(0,\ldots,0)=w_k$, 
$1\le k\le n-1$. It follows from~\eqref{c-ex1} that if $\varepsilon_1\neq 0$  
then there exists $j\in \{1,\ldots,n-1\}$ such that 
\begin{equation}\label{realpart}
|\arg \varepsilon_1-\arg w_j|\le \frac{\pi}{n-1},\text{ so that } 
\Re\left(\frac{\varepsilon_1}{w_j}\right)\ge 
\frac{|\varepsilon_1|}{d(p)}\cos\!\left(\frac{\pi}{n-1}\right).
\end{equation}
Now using Lemma~\ref{invfcn} with $v=(0,z_2,\ldots,z_n)$ and 
the computations in the proof of Theorem~\ref{thm15} together with the 
assumption
$|\varepsilon_i|\le |\varepsilon_1|^{1+\kappa}$, $2\le i\le n$, we get 
\begin{equation}\label{expansion}
\begin{split}
\left|\omega_j(\varepsilon_1,\ldots,\varepsilon_n)-\varepsilon_1\right|
&=\left|w_j-\varepsilon_1
+\sum_{i=1}^{n}\frac{\partial \omega_j}{\partial \zeta_i}\bigg|_{v}
\varepsilon_i+\mathcal{O}\left(|\varepsilon_1|^2\right)\right|\\
&=|w_{j}|-\frac{n-1}{n}\Re\left(\frac{\varepsilon_1}{w_j}\right)+
\mathcal{O}\left(|\varepsilon_1|^{1+\kappa'}\right),
\end{split}
\end{equation}
where $\kappa'=\min(1,\kappa)$. It is not difficult to see that if $n\ge 4$ 
then $d(p)=\left(\!\frac{1}{n}\!\right)^{\!\frac{1}{n-1}}<1-\frac{1}{n}$. 
From~\eqref{realpart}--\eqref{expansion} and the latter inequality it follows 
that if $n\ge 4$ then for all small $|\varepsilon_1|$ one has
\begin{equation*}
\begin{split}
\left|\omega_j(\varepsilon_1,\ldots,\varepsilon_n)-\varepsilon_1\right|
&\le |w_{j}|-\left[\frac{n-1}{nd(p)}\cos\!\left(\frac{\pi}{n-1}\right)\right]
|\varepsilon_1|+\mathcal{O}\left(|\varepsilon_1|^{1+\kappa'}\right)\\
&\le d(p)-\left[\cos\!\left(\frac{\pi}{n-1}\right)\right]|\varepsilon_1|,
\end{split}
\end{equation*}
which proves the proposition.
\end{proof}

As already noted in Remark~\ref{remmax2}, if $p\in S_n$ and $\ze\in Z(p)$ is 
such that $|\ze|=1$ then $|p|_{\ze}\le 1$. It is therefore natural to ask 
whether the following holds:

\begin{problem}\label{pb2}
Let $p(z)=(z-z_1)q(z)\in S_n$, where $|z_1|<1$ and $q\in S_{n-1}$. Does there 
exist $\ze\in\bC$ with $|\ze|=1$ such that $|p|_{z_1}\le |r|_{\ze}$, where 
$r(z)=(z-\ze)q(z)$?
\end{problem}

In other words, is it always possible to increase the Sendov distance for a 
given zero of $p$ in the unit disk by ``pushing'' this zero to the boundary 
while keeping the other zeros of $p$ fixed? The $0$-maximal polynomials of the 
form $z^n-z$ show that such a property cannot hold in general:

\begin{proposition}\label{v8prop1}
If $n\ge 5$ there exist $z_1\in D$ and $q\in S_{n-1}$ such that for any 
$\ze\in\bC$ with $|\ze|=1$ one has 
$|r|_{\ze}<|p|_{z_1}$, where $p(z)=(z-z_1)q(z)$ and $r(z)=(z-\ze)q(z)$.
\end{proposition}

\begin{proof}
Let $\ze\in \mathbf{C}$ be such that $|\ze|=1$ and set $z_1=0$ and 
$q(z)=z^{n-1}-1$, so that $p(z)=z^n-z$ and 
$r(z)=(z-\ze)(z^{n-1}-1)$. Then for any $\omega\in Z(p')$ one has
$|z_1-\omega|=\left(\frac{1}{n}\right)^{\frac{1}{n-1}}$, so that
\begin{equation}\label{v8eq2}
|p|_{z_1}=\min_{\omega\in Z(p')}|z_1-\omega|
=\left(\frac{1}{n}\right)^{\frac{1}{n-1}}.
\end{equation}
Since $|\ze|=1$ and $Z(q)=\{e^{\frac{2k\pi i}{n-1}}\,:\,0\le k\le n-2\}$ one 
gets from either the Schur-Szeg\"o composition theorem 
(\cite[Theorem 3.4.1d]{RS}) or the Grace-Heawood theorem 
(\cite[Theorem 4.3.1]{RS}) that
\begin{equation}\label{v8eq3}
|r|_{\ze}=\min_{\omega\in Z(r')}|\ze-\omega|\le \frac{\min_{\alpha\in Z(q)}
|\ze-\alpha|}{2\sin(\pi/n)}\le \frac{\sin\left(\pi/2(n-1)\right)}{\sin(\pi/n)}.
\end{equation}
From~\eqref{v8eq2} and~\eqref{v8eq3} we see that in order to prove the 
proposition it is enough to show that
\begin{equation}\label{v8eq4}
\frac{\sin\left(\pi/2(n-1)\right)}{\sin(\pi/n)}
< \left(\frac{1}{n}\right)^{\frac{1}{n-1}}\text{ for }n\ge 5.
\end{equation}
Numerical checking shows that~\eqref{v8eq4} is true for $n=5,6,7$ or 8, and so 
we may assume that $n\ge 9$. Since 
$0<\dfrac{\pi}{2(n-1)}<\dfrac{\pi}{n}<\dfrac{\pi}{4}$ for $n\ge 5$ and 
$x\mapsto \dfrac{\sin x}{\sqrt{x}}$ is an 
increasing function on $\left(0,\dfrac{\pi}{4}\right)$ we get that
\begin{equation}\label{v8eq5}
\frac{\sin\left(\pi/2(n-1)\right)}{\sin(\pi/n)}<\sqrt{\frac{n}{2(n-1)}}
\,\text{ for }n\ge 5.
\end{equation}
The sequence $\sqrt{2\left(1-\frac{1}{n}\right)}-n^{\frac{1}{n-1}}$ is clearly 
increasing for $n\ge 2$, so that 
$$\sqrt{2\left(1-\frac{1}{n}\right)}-n^{\frac{1}{n-1}}\ge 
\frac{4}{3}-3^{1/4}>0$$ 
whenever $n\ge 9$. This implies that the right-hand side of~\eqref{v8eq5} is 
always less than $\left(\frac{1}{n}\right)^{\frac{1}{n-1}}$ if $n\ge 9$, which 
proves~\eqref{v8eq4}. 
\end{proof}

\section{Second order variations and locally maximal polynomials}\label{s-3}

Theorem~\ref{thm15} and Lemma~\ref{countlm3} show that $z^3+z$ is a linearly 
inextensible polynomial which is not a local maximum for $d$ in $S_n$. On the 
other hand, all the properties of the polynomial 
$p(z)=z^n+z$ that we discussed so far seem to suggest that if $n\ge 4$ then 
$p$ and its rotations could in fact be locally maximal for Sendov's 
conjecture in the whole class $S_n$. More precisely:

\begin{enumerate}
\item[(i)] $p$ is locally maximal for the function $d|_{S(n,0)}$ 
(Proposition~\ref{locmax0}); 
\item[(ii)] $p$ is linearly inextensible (Theorem~\ref{thm15});
\item[(iii)] $p$ is locally maximal for many types of variations 
(Proposition~\ref{propvar4});
\item[(iv)] the zeros and critical points of $p$ are symmetrically 
distributed and satisfy Conjecture~\ref{miller} (Remark~\ref{remmax1}).
\end{enumerate}

We shall now use a second order 
variational method to show that -- contrary to what one might expect from 
properties (i)--(iv) -- the polynomial $p$ is {\em not} locally maximal for 
Sendov's conjecture. 
Theorems~\ref{secord4} and~\ref{secord5} below show that if $n=4$ or $n=5$ 
then 
$p$ is a kind of inflection point for the function $d$. We conjecture that the 
same is actually true for all degrees (Remark~\ref{gencase}) and we also 
discuss the problem of finding the local maxima for $d$ in 
$S_n$ (Problem~\ref{pb3}).

\begin{notation}\label{not-4}
Let $a\in [0,1]$ and set
\begin{equation}\label{var1}
\begin{split}
&r=\left|z^4+z\right|_0=\sqrt[3]{\frac{1}{4}},
\quad \alpha_1=\frac{3\sqrt{3}r}{2-3r},\quad
\alpha_2=-\frac{\sqrt{3}\left[(3r+2)^2+4\right]}{2(3r-2)^2},\\
&\zeta(a)=\exp{\left[i\!\left(\frac{\pi}{3}+\alpha_1 a
+\alpha_2 a^2\right)\right]}.
\end{split}
\end{equation}
\end{notation}

We use the quantities in~\eqref{var1} in order to construct certain continuous 
deformations of the polynomial $z^4+z$. For sufficiently small $a>0$ these  
perturbations give rise to a one-parameter 
family of polynomials in $S_4$ with the following property:

\begin{theorem}\label{secord4}
Let $\zeta(a)$ be as in~\eqref{var1} and define the polynomial 
\begin{equation}\label{var2}
p_a(z):=(z-a)(z+1)\left[z-\zeta(a)\right]\!
\left[z-\overline{\zeta(a)}\right]\in S_4.
\end{equation}
Then for all sufficiently small $a>0$ one has 
$\Delta(p_a,z^4+z)=\mathcal{O}(a)$ and 
$$d(p_a)=|p_a|_a=r+Ca^2+\mathcal{O}(a^3),\,\text{ where }\,
C=\frac{3}{4r(2-3r)}\approx 10.81154938.$$
In particular, the polynomial $z^4+z$ is not locally maximal for 
Sendov's conjecture.
\end{theorem}

\begin{proof}
Let us first describe the steps we took in order to arrive at the 
quantities in~\eqref{var1} and the one-parameter family of polynomials 
defined in~\eqref{var2}. We start by introducing six real parameters 
which we denote by $x_i$, $y_i$, $1\le i\le 3$, and we define the auxiliary 
polynomials
\begin{eqnarray*}
P_a(z;x_1,x_2,x_3)\!&=&\!(z-a)\left[z-z_1(a;x_3)\right]
\left[z-z_2(a;x_1,x_2)\right]
\!\left[z-\overline{z_2(a;x_1,x_2)}\right]\\
&=&\!z^4+\sum_{j=0}^{3}b_j(a;x_1,x_2,x_3)z^j,\\
Q_a(z;y_1,y_2,y_3)\!&=&\!4\!\left[z-\omega_1(a;y_3)\right]
\left[z-\omega_2(a;y_1,y_2,y_3)\right]
\!\left[z-\overline{\omega_2(a;y_1,y_2,y_3)}\right]\\
&=&\!4z^3+\sum_{k=0}^{2}c_k(a;y_1,y_2,y_3)z^k,
\end{eqnarray*}
where
\begin{equation}\label{omegas}
\begin{split}
&z_1(a;x_3)=-1+x_3 a^2,\quad z_2(a;x_1,x_2)=\exp{\left[i\!\left(\frac{\pi}{3}
+x_1 a+x_2 a^2\right)\right]},\\
&\omega_1(a;y_3)=a-\left(r+y_3 a^2\right),\\
&\omega_2(a;y_1,y_2,y_3)=
a+\left(r+y_3 a^2\right)
\exp{\left[i\!\left(\frac{\pi}{3}+y_1 a+y_2 a^2\right)\right]}.
\end{split}
\end{equation}
The idea is now to investigate whether the parameters 
$x_i$ and $y_i$, $1\le i\le 3$, may be chosen so that the following 
conditions are satisfied:
\begin{itemize}
\item[(i)] $x_3\ge 0$ and $y_3>0$;
\item[(ii)] for all sufficiently small $a>0$ one has
$$\max_{|z|\le 1}\big|P'_a(z;x_1,x_2,x_3)-Q_a(z;y_1,y_2,y_3)\big|
=\mathcal{O}(a^3),$$
\end{itemize}
where $P'_a(z;x_1,x_2,x_3)$ denotes the derivative of $P_a(z;x_1,x_2,x_3)$ 
with respect to $z$. To do this, we expand the coefficients of 
$P'_a(z;x_1,x_2,x_3)$ and $Q_a(z;y_1,y_2,y_3)$ into their Maclaurin series 
so as to get second order approximations of these coefficients (with an error 
of $\mathcal{O}(a^3)$). Let $\tilde{b}_j(a;x_1,x_2,x_3)$, $1\le j\le 3$, and 
$\tilde{c}_k(a;y_1,y_2,y_3)$, $0\le k\le 2$, be the resulting second degree 
Maclaurin polynomials in the variable $a$ for the coefficients 
$b_j(a;x_1,x_2,x_3)$, $1\le j\le 3$, and $c_k(a;y_1,y_2,y_3)$, $0\le k\le 2$, 
respectively. Then we may write
\begin{multline}
(m+1)\tilde{b}_{m+1}(a;x_1,x_2,x_3)-\tilde{c}_{m}(a;y_1,y_2,y_3)\\
=\sum_{n=0}^{2}d_{mn}(x_1,x_2,x_3,y_1,y_2,y_3)a^n,\quad 0\le m\le 2,
\end{multline}
where $d_{mn}(x_1,x_2,x_3,y_1,y_2,y_3)$, $0\le m,n\le 2$, are 
real polynomials in the variables $x_1,x_2,x_3,y_1,y_2,y_3$.
Clearly, any solution $(x_1,x_2,x_3,y_1,y_2,y_3)\in \mathbf{R}^6$ to the 
system of polynomial equations
\begin{equation}\label{syst96}
d_{mn}(x_1,x_2,x_3,y_1,y_2,y_3)=0,\quad 0\le m,n\le 2,
\end{equation}
that satisfies $x_3\ge 0$ and $y_3>0$ will also satisfy conditions (i) and 
(ii) above. It turns out that if we let $x_3=0$ then system~\eqref{syst96} 
may be reduced to a system of five linear equations in the variables 
$x_1,x_2,y_1,y_2,y_3$ which admits a unique solution. We arrive in this way 
at the following solution to system~\eqref{syst96}:
\begin{equation}\label{solsyst96}
\begin{split}
&x_1=\alpha_1,\quad x_2=\alpha_2,\quad x_3=0,\quad 
y_1=\frac{3\sqrt{3}}{2r(2-3r)},\\
&y_2=-\frac{3\sqrt{3}(12r^2+8r+3)}{8r^2(3r-2)^2},
\quad y_3=C:=\frac{3}{4r(2-3r)},
\end{split}
\end{equation}
where $r$, $\alpha_1$ and $\alpha_2$ are as in~\eqref{var1}. Henceforth we 
assume that the values of the parameters $x_1,x_2,x_3,y_1,y_2,y_3$ are those 
listed in~\eqref{solsyst96}. In particular, this implies that 
$P_a(z;x_1,x_2,x_3)$ is the same as the polynomial $p_a(z)$ given 
in~\eqref{var2}. Note that by~\eqref{var1} one has 
$\Delta(p_a,z^4+z)=\mathcal{O}(a)$ for all small positive $a$. To simplify 
the notations, the zeros of the polynomial $Q_a(z;y_1,y_2,y_3)$ will be 
denoted by $\omega_1(a)$, $\omega_2(a)$ and $\overline{\omega_2(a)}$, 
respectively. Then~\eqref{omegas} and~\eqref{solsyst96} imply that
\begin{equation}\label{concl}
|\omega_1(a)-a|=|\omega_2(a)-a|=r+Ca^2.
\end{equation}

It is now practically clear that the desired conclusion should follow from 
the Newton-Raphson algorithm. We check this in a rigorous way by using 
the following simple observation.

\begin{lemma}\label{NR}
Let $R$ be a complex polynomial of degree $d\ge 2$. If $w\in \mathbf{C}$ is 
such that $R'(w)\neq 0$ then there exists a zero $z$ of $R$ such that
$$|w-z|\le d\left|\frac{R(w)}{R'(w)}\right|.$$
\end{lemma}

\begin{proof}
If $R(w)=0$ there is nothing to prove. Otherwise, we let 
$z_1,\ldots,z_d$ denote the zeros of $R$. Then
$$\sum_{i=1}^{d}\frac{1}{|w-z_i|}\ge \left|\sum_{i=1}^{d}\frac{1}{w-z_i}\right|
=\left|\frac{R'(w)}{R(w)}\right|,$$
which proves the lemma.
\end{proof}

Let $w_1(a)$, $w_2(a)$ and $w_3(a)$ denote the critical points of the 
polynomial $p_a(z)$. Since $p_a(z)\rightarrow z^4+z$ as $a\rightarrow 0$, 
we may label these critical points so that $w_1(a)\in \mathbf{R}$, 
$\Im (w_2(a))>0$ and $w_3(a)=\overline{w_2(a)}$ if $a$ is positive and 
sufficiently small. A straightforward computation shows that
\begin{equation*}
p_{a}''(\omega_1(a))=\frac{3}{r}+\mathcal{O}(a)\,\text{ and }\,
p_{a}''(\omega_2(a))=\overline{p_{a}''\left(\overline{\omega_2(a)}\right)}
=\frac{3e^{\frac{2\pi i}{3}}}{r}+\mathcal{O}(a),
\end{equation*}
which combined with Lemma~\ref{NR} and condition (ii) implies that
$$w_1(a)=\omega_1(a)+\mathcal{O}(a^3)\,\text{ and }\,w_2(a)=\overline{w_3(a)}
=\omega_2(a)+\mathcal{O}(a^3).$$
The desired result is now a consequence of~\eqref{concl}. 
\end{proof}

As one may expect, the second order variational method that we used in the 
proof of Theorem~\ref{secord4} works even in more general cases. However, for 
higher degrees the procedure described above requires an increasingly large 
amount of computations. Tedious as they may be when done by hand, for small 
degrees these computations become relatively easy if one uses for instance a 
Maple computer program. Such a program has considerably simplified 
our task in the course of proving Theorem~\ref{secord4} and it also allowed 
us to obtain a similar result for the polynomial $z^5+z$. 

\begin{notation}\label{not-5}
Let $a\in [0,1]$ and set 
\begin{equation*}
\begin{split}
&s=\left|z^5+z\right|_0=\sqrt[4]{\frac{1}{5}},\quad
\beta=\frac{2\sqrt{2}s^2}{1-2s^2},\quad
\gamma=\frac{4\sqrt{2}}{5s(1-2s^2)},\\
&\delta=\frac{60s^4-19}{50s^2(2s^2-1)^2},\quad 
K=\frac{2}{5s(1-2s^2)}\approx 5.665658792,\\
&\eta(a)=\exp{\left[i\!\left(\frac{\pi}{4}+\beta a\right)\right]},\quad 
\chi_1(a)=a+\left(s+Ka^2\right)
\exp{\left[i\!\left(\frac{\pi}{4}+\gamma a+\delta a^2\right)\right]},\\
&\chi_2(a)=a+\left(s+Ka^2\right)
\exp{\left[i\!\left(\frac{3\pi}{4}+\gamma a-\delta a^2\right)\right]}.
\end{split}
\end{equation*}
We use these quantities in order to define the following one-parameter 
families of polynomials:
\begin{eqnarray*}
&&q_a(z):=(z-a)(z-\eta(a))(z-i\eta(a))\!\left[z-\overline{\eta(a)}\right]
\!\left[z+i\overline{\eta(a)}\right]\in S_5,\\
&&s_a(z):=5(z-\chi_1(a))(z-\chi_2(a))\!\left[z-\overline{\chi_1(a)}\right]
\!\left[z-\overline{\chi_2(a)}\right].
\end{eqnarray*}
\end{notation}
We can now formulate an analogue of Theorem~\ref{secord4} for the polynomial 
$z^5+z$. Its proof is similar to that of Theorem~\ref{secord4} and is 
therefore omitted.

\begin{theorem}\label{secord5}
For all sufficiently small $a>0$ one has
\begin{equation*}
\Delta(q_a,z^5+z)=\mathcal{O}(a)\,\text{ and }\,
\Delta\!\left(\frac{q'_a}{5},\frac{s_a}{5}\right)=\mathcal{O}(a^3),
\end{equation*}
so that $d(q_a)=|q_a|_a=s+Ka^2+\mathcal{O}(a^3)$. In particular, 
the polynomial $z^5+z$ is not locally maximal for Sendov's conjecture.
\end{theorem}

\begin{remark}\label{gencase}
In view of Lemma~\ref{countlm3} and Theorems~\ref{secord4}--\ref{secord5} it 
seems reasonable to conjecture that polynomials of 
the form $z^n+z$ or, more generally, the $0$-maximal polynomials 
of degree $n\ge 2$ given in Theorem~\ref{class0max} are not local maxima for 
the function $d$. One could for instance try to find an algorithmic proof of 
this conjecture based on the 
second order variational method that we just described. 
\end{remark}

\begin{remark}\label{limit}
Theorems~\ref{secord4}--\ref{secord5} and the second order variational methods
used in \cite{M3}, \cite{M4} and \cite{VZ} clearly expose the limitations 
of what can be achieved through first order expansions and show the 
necessity of using
higher order approximations when dealing with problems such as Sendov's 
conjecture. 
\end{remark}

The polynomials in Conjecture~\ref{send2} -- i.e., the conjectured extremal 
polynomials for Sendov's conjecture -- have long been the only 
known examples of local maxima for $d$ in $S_n$ (cf.~\cite{M3} and \cite{VZ}). 
Moreover, Sendov's conjecture and related problems have been studied almost 
exclusively by means of analytical or variational methods and all the results 
obtained until very recently seemed to suggest a negative answer to the 
following question:

\begin{problem}\label{pb3}
Are there any local maxima for $d$ in $S_n$ other than $z^n+e^{i\theta}$,
$\theta\in\bR$?
\end{problem}

As pointed out in \cite{B}, the fact that $d$ fails to be a (logarithmically) 
plurisubharmonic function in the polydisk $\bar{D}^n$ accounts 
for many of the difficulties in answering Problem~\ref{pb3}. In a remarkable 
recent paper \cite{M5} Miller showed that the answer to this question 
is in fact affirmative by explicitly constructing local maxima of degree 
$n\in \{8, 12, 14, 20, 26\}$ which are ``unexpected'' in the sense that they 
are not of the form $z^n+a$ with $|a|=1$. In particular, Miller's result 
shows that local variational methods alone (of {\em any} order) are not 
enough for dealing with Conjecture~\ref{send2} in full generality.
It is therefore important to consider alternative 
methods for studying Sendov's conjecture and its strong version 
(Conjecture~\ref{send2}). An inductive 
approach to Conjectures~\ref{send1} and~\ref{send2} based on apolarity and 
coincidence theorems of Grace-Walsh type was discussed in~\cite{B}. More 
recently, operator theoretical interpretations and extensions of Sendov's 
conjecture and related questions were proposed in~\cite{B1,BP1,BP2}. In 
{\em loc.~cit.~}it was also shown that Conjectures~\ref{send1} 
and~\ref{send2} may in fact be viewed as part of the more general problem of 
describing the relationships between the spectra of normal matrices and the 
spectra of their degeneracy one principal submatrices. These and similar 
questions form the basis of an ongoing joint study with R.~Pereira and 
S.~Shimorin and will be addressed in a forthcoming paper.


\begin{thebibliography}{99}

\bibitem{BRS}
Bojanov, B., Rahman, Q.~I., Szynal, J., {\em On a conjecture about the 
critical points of a polynomial}, in ``Delay Equations, Approximation and 
Application'' (G.~Meinardus, G.~N\"urnberger, Eds.), 83--93, 
Internat. Series of 
Numer. Math., Vol. 74, Birkh\"auser Verlag, Basel, 1985. 

\bibitem{B}
Borcea, J., {\em Two approaches to Sendov's conjecture}, 
Arch.~Math.~{\bf 71} (1998), 46--54.

\bibitem{B1}
Borcea, J., {\em Maximal and inextensible polynomials and the geometry of the 
spectra of normal operators}, arXiv:math.CV/0309233.

\bibitem{BP1}
J.~Borcea, M.~Pe\~na, {\em Equilibrium points of logarithmic potentials 
induced by positive charge distributions. I. Generalized de Bruijn-Springer 
relations}, Trans. Amer. Math. Soc., to appear; arXiv:math.CV/0601519.

\bibitem{BP2}
J.~Borcea, M.~Pe\~na, {\em Equilibrium points of logarithmic potentials 
induced by positive charge distributions. II. A conjectural Hausdorff 
geometric symphony}, submitted; preprint (2005).

\bibitem{BX}
Brown, J.~E, Xiang, G., {\em Proof of the Sendov conjecture for polynomials of 
degree at most eight}, J.~Math.~Anal.~Appl.~{\bf 232} (1999), 272--292.

\bibitem{M1}
Miller, M.~J., {\em Maximal polynomials and the Ilieff-Sendov conjecture}, 
Trans. Amer. Math. Soc. {\bf 321} (1990), 285--303.

\bibitem{M2}
Miller, M.~J., {\em Continuous independence and the Ilieff-Sendov 
conjecture}, Proc. Amer. Math. Soc. {\bf 115} (1992), 79--83.

\bibitem{M3}
Miller, M.~J., {\em On Sendov's conjecture for roots near the unit circle}, 
J. Math. Anal. Appl. {\bf 175} (1993), 632--629.

\bibitem{M4}
Miller, M.~J., {\em A quadratic approximation to the Sendov radius near the 
unit circle}, Trans. Amer. Math. Soc. {\bf 357} (2005), 851--873.  

\bibitem{M5}
Miller, M.~J., {\em Unexpected local extrema for the Sendov conjecture}, 
arXiv:math.CV/0505424.

\bibitem{PR}
Phelps, D., Rodriguez, R.~S., {\em Some properties of extremal polynomials 
for the Ilieff conjecture}, K\^odai Math.~Sem.~Rep.~{\bf 24} (1972), 172--174.

\bibitem{RS}
Rahman, Q.~I., Schmeisser, G., Analytic theory of polynomials, London 
Math. Soc. Monogr. (N.~S.) Vol. 26, Oxford Univ.~Press, New York, NY, 2002.

\bibitem{Ro}
Rockafellar, R.~T., Convex analysis, Princeton Univ.~Press, Princeton, NJ, 
1997.

\bibitem{Ru}
Rubinstein, Z., {\em On a problem of Ilyeff}, Pacific J.~Math.~{\bf 26} 
(1968), 159--161.

\bibitem{Sc}
Schmeisser, G., {\em The conjectures of Sendov and Smale}, in ``Approximation 
Theory: A volume dedicated to Blagovest Sendov'' (B.~Bojanov, Ed.), 353--369. 
Darba, Sofia, 2002. 

\bibitem{SSz}
Schmieder, G., Szynal, J., {\em On the distribution of the derivative zeros 
of a complex polynomial}, Complex Var.~Theory Appl.~{\bf 47} (2002), 239--241.

\bibitem{Se}
Sendov, B., {\em Hausdorff geometry of polynomials}, East J.~Approx.~{\bf 7} 
(2001), 123--178.

\bibitem{Sh}
Sheil-Small, T., Complex Polynomials, Cambridge Studies 
in Adv.~Math.~vol.~{\bf 75}, Cambridge Univ.~Press, Cambridge, UK, 2002.

\bibitem{VZ}
V\^aj\^aitu, V., Zaharescu, A., {\em Ilyeff's conjecture on a corona}, 
Bull.~London Math.~Soc.~{\bf 25} (1993), 49--54.

\end{thebibliography}
\end{document}